\documentclass[12pt]{article}	

\long\def\RTransNotes#1\EndRTransNotes{
}
\def\Title{Constructing new Calabi--Yau $3$-folds and their 
\\[9pt]		mirrors via conifold transitions}

\def\Abstract{ ~ ~
	We construct a surprisingly large class of new Calabi--Yau $3$-folds 
$X$ with small Picard numbers and propose a construction of their 
mirrors $X^*$ using smoothings of toric hypersurfaces with conifold 
singularities. These new examples are related to the previously known 
ones via conifold transitions. Our results generalize the 
mirror construction for Calabi--Yau complete intersections in 
Grassmannians and flag manifolds via toric degenerations. 
There exist exactly $198849$ reflexive $4$-polytopes whose 2-faces are 
only triangles or parallelograms of minimal volume. Every such polytope  
gives rise to a family of Calabi-Yau hypersurfaces with at worst conifold 
singularities. Using a criterion of Namikawa we found $30241$  
reflexive $4$-polytopes such that the corresponding Calabi-Yau 
hypersurfaces are smoothable by a flat deformation. In particular, 
we found $210$ reflexive $4$-polytopes defining $68$ topologically 
different Calabi--Yau $3$-folds with $h_{11}=1$. We explain the mirror 
construction and compute several new Picard--Fuchs operators for the
respective $1$-parameter families of mirror Calabi-Yau $3$-folds.
}

\usepackage{amssymb,epsfig,epic}	\long\def\del#1\enddel{} 

\include{pst-plot}		\unitlength=3pt		\psset{unit=3pt}
\def\BP{\begin{picture}}	\def\EP{\end{picture}}	\thicklines
\def\BI{\begin{itemize}}	\def\EI{\end{itemize}}
\def\BE{\begin{equation}}	\def\EE{\end{equation}}
\def\BEA{\begin{eqnarray}}   	\def\EEA{\end{eqnarray}} 

\def\fnote#1#2{\begingroup\def\thefootnote{#1}\footnote{#2}
                \addtocounter{footnote}{-1}\endgroup}
\def\ifundefined#1{\expandafter\ifx\csname#1\endcsname\relax}

\def\IZ{{\mathbb Z}}	\def\IR{{\mathbb R}} 	\def\IP{{\mathbb P}}	
\def\IC{{\mathbb C}}
\def\IT{{\mathbb T}}				\let\ex=\times
\def\IQ{{\mathbb Q}}				\let\PP=\varpi

\let\th=\theta	\let\D=\Delta	\let\L=\Lambda	\let\S=\Sigma

\begin{document}

\begin{center}
\vspace*{12mm}{\Large\bf	\Title			}
\\[12mm]
\large 				Victor BATYREV\fnote{*}
{E-mail: \it victor.batyrev@uni-tuebingen.de}	\normalsize
\\[4mm]	Mathematisches Institut, Universit\"at T\"ubingen
\\[2mm]	Auf der Morgenstelle 10, 72076 T\"ubingen, GERMANY
\\[4mm]	and
\\[4mm]				Maximilian KREUZER\fnote{\#}
{E-mail: \it maximilian.kreuzer@tuwien.ac.at}
\\[4mm]	Institute for Theoretical Physics, Vienna Univ. of Technology
\\[2mm]	Wiedner Hauptstr. 8-10, A-1040 Vienna, AUSTRIA
\\[3cm]				\bf Abstract
\\[9mm]\parbox{150truemm}{	\rm \Abstract	}
\end{center}			\thispagestyle{empty}

\newpage		\baselineskip=16pt

\section{Introduction}

Toric geometry provides simple and efficient combinatorial tools 
\cite{BatNill} for the construction of large classes of 
Calabi--Yau manifolds from generic hypersurfaces \cite{Bat} 
and complete intersections \cite{BBnef,BBsth} in toric varieties. An 
important additional benefit of this construction  is its  invariance
under mirror symmetry. In particular it  enables the
computation of quantum cohomology and instanton numbers using 
generalized hypergeometric functions \cite{BatQC,CoKa}.
Calabi--Yau $3$-folds obtained from hypersurfaces in 
$4$-dimensional toric varieties 
have been enumerated 
completely \cite{c4d,pwf}. Some  large lists of Calabi--Yau $3$-folds 
obtained from complete intersections
 have been compiled and analyzed in \cite{wpci,kkrs}. 
Their  fibration structures
\cite{Tfib,pwf}
and torsion in cohomology \cite{IC} are  of particular interest for 
applications to string theory \cite{kkrs,WSI}.
Thus, toric constructions provide by far the largest number of known
examples, but they are, nevertheless, quite special in the zoo of all 
Calabi-Yau $3$-folds  about which little is known. Even  finiteness of 
topological types of Calabi-Yau $3$-folds remains still an open question.
 
According to Miles Reid \cite{Reid}
it is expected  that an appropriate partial  compactification of  
the  moduli space of all Calabi--Yau $3$-folds, which  allows Calabi-Yau 
varieties  with   mild singularities, will be connected. 
Using this idea,  we can try to get new examples of Calabi-Yau 
$3$-folds  by studying singular 
limits of Calabi-Yau $3$-folds obtained by toric
methods.

In the present work we focus on Calabi-Yau $3$-folds $\widehat{X}_f$ 
obtained from generic hypersurfaces $\overline{X}_f$ 
in $4$-dimensional Gorenstein toric Fano varieties $\IP_\D$ correspondig to 
$4$-dimensional reflexive polytopes $\D$ 
because the  complete list of these polytopes  
is known \cite{c4d,pwf}. We classify all hypersurfaces  $\overline{X}_f$ 
with at most 
conifold singularities coming from the singularities of the ambient 
Gorenstein toric Fano variety $\IP_\D$. Standard toric methods 
allow to resolve  the conifold 
singularities of  $\overline{X}_f$ by a toric resolution  of 
 $\IP_\D$ and to obtain a smooth Calabi-Yau $3$-fold $\widehat{X}_f$. 
However, in this paper we are interested in smoothing  $\overline{X}_f$ to 
a Calabi-Yau $3$-fold $Y$ by a
flat deformation. Thus the two Calabi-Yau $3$-folds 
 $\widehat{X}_f$ and $Y$ are connected by a so-called {\em conifold 
 transition}. Moreover, a similar 
conifold transition exists for the mirrors 
 $(\widehat{X}_f)^*$ and $Y^*$.

The above  approach was  motivated by the  
previous work \cite{Grass,Flag} which shows that 
all Grassmannians  and Flag manifolds allow flat degenerations
to Gorenstein toric Fano varieties $\IP_\D$ having at worst conifold 
singularities in codimension $3$. Therefore smooth $3$-dimensional 
Calabi-Yau complete intersections $Y$ in  
Grassmannians  and Flag manifolds 
can be regarded as smoothings of generic Calabi-Yau 
complete intersections $\overline{X}$ in the corresponding 
Gorenstein toric Fano variety $\IP_\D$.  
It was observed in  \cite{Grass,Flag} that the mirrors $Y^*$ of 
$Y$ are   obtained by  specializations of the complex structure 
of the mirrors  $\widehat{X}^*$ coming from the  already known 
toric construction.

In section 2 we explain the construction and present our results by describing
the lists of polytopes and Hodge data whose details are available on the
internet. 
In section 3 we discuss the conifold transition in the mirror
family (which is related to the transitions studied in \cite{Rolling}).
This enables the construction of non-toric mirror pairs and the computation 
of quantum cohomologies.  
In section 4 we focus on 1-parameter models for which we compute the
topological data and also initiate the study of the mirror map by direct 
evaluation of the principal period, which allows us to find a number of 
new Picard-Fuchs operators.
We conclude with a discussion of open problems, generalizations,
and work to be done.

\section{Reflexive polytopes and conifold transitions}

Let $M\cong \IZ^d$ and $N=\hbox{Hom}(M,\IZ)$ be a dual pair of lattices 
together with the canonical pairing $\langle *,*\rangle\; : \; 
N \times M \to \IZ$ 
and let 
$M_\IR=M\otimes\IR$, $N_\IR=N\otimes\IR$ be 
their  real extensions. It is known \cite{Bat} that   
the generic families of  
Calabi--Yau hypersurfaces $\overline{X}_f$ in $4$-dimensional 
Gorenstein toric Fano varieties $\IP_\D$ 
and their mirrors  $\overline{X}_g \subset \IP_{\D^\circ}$
are in one-to-one correspondence
to the polar pairs $\D\subset M_\IR$, $\D^\circ\subset N_\IR$ 
of reflexive 4-polytopes.
By definition reflexivity of $\D$ and $\D^\circ$ means that 
\BE							\label{polar}
	\D^\circ=\{y\in N_\IR~:~\langle y,x\rangle\ge-1~\forall x\in\D\}
\EE
and that both $\D$ and  $\D^\circ$ are lattice polytopes, i.e., 
all vertices of $\D$ (resp. $\D^\circ$) are elements 
of $M$ (resp. $N$). 

The hypersurfaces $\overline{X}_f \subset \IP_\D$ and   
$\overline{X}_g \subset \IP_{\D^\circ}$ 
are  the closures 
of the affine hypersurfaces $X_f$ and $X_g$ defined by generic Laurent 
polynomials 
\BE                            \label{CYeq}
 f:= \sum_{m \in \D \cap M} a_m t^m \;\quad\mbox{\rm and} \;\quad
 g:= \sum_{n \in \D^{\circ} \cap N} b_n t^n.   
\EE
Denote by $\Sigma$ (resp. by $\Sigma^{\circ}$)  the fan of cones 
over simplices in $\partial \Delta$ (resp. in  $\partial \Delta^{\circ} $)
in a maximal coherent triangulation 
of $\Delta$ (resp. $\D^{\circ})$  \cite{GKZ}. Then the projective 
toric variety $\IP_{\Sigma^{\circ}}$ is a maximal partial projective 
crepant (MPPC) desingularization of $\IP_\D$. Similarly, 
 $\IP_{\Sigma}$ is a maximal partial projective 
crepant desingularization of $\IP_{\D^{\circ}}$. Denote by 
$\widehat{X}_f \subset \IP_{\Sigma^{\circ}}$ and $\widehat{X_g} \subset 
\IP_{\Sigma}$ the closures of 
$X_f$ and  $X_g$ in   $\IP_{\Sigma^{\circ}}$ and  $\IP_{\Sigma}$,
respectively. Then $\widehat{X}_f$ and $\widehat{X_g}$ are
 smooth Calabi--Yau 3-folds   \cite{Bat} and for the Hodge numbers 
$h^{1,1}$ and $h^{2,1}$ 
one has 
\BE                            \label{MPC}
     h^{1,1}(\widehat{X}_f)=h^{2,1}(\widehat{X}_g) =
                l(\D^\circ)-5
   	-\!\!\!\!\sum_{\mathrm{codim}(\th^\circ)=1}l^*(\th^\circ)
                +\!\!\!\!
       \sum_{\mathrm{codim}(\th^\circ)=2}l^*(\th^\circ)l^*(\th), \EE
\BE h^{1,1}(\widehat{X}_g)=h^{2,1}(\widehat{X}_f) =
                l(\D)-5
   	-\!\!\!\!\sum_{\mathrm{codim}(\th)=1}l^*(\th)
                +\!\!\!\!
       \sum_{\mathrm{codim}(\th)=2}l^*(\th)l^*(\th^\circ), 
\EE
where $l(\D)$ denotes the number of lattice points of $\D$ and $l^*(\th)$
denotes the number of lattice points in the relative interior of $\th$. 
The faces $\th\subset\D$ and
$\th^\circ\subset\D^\circ$ denote polar sets of points saturating the 
inequality in eq. (\ref{polar}), so that $\dim(\th)+\dim(\th^\circ)=3$.

The smoothness of generic hypersurfaces $\widehat{X}_f$ and $\widehat{X}_g$ 
follows from the fact that the singularities of the MPPC resolutions 
$\IP_{\Sigma^{\circ}}$ and $\IP_{\Sigma}$ 
have codimention at least 4, i.e., singular  points
in the  4-dimensional ambient spaces  $\IP_{\Sigma^{\circ}}$ and 
$\IP_\S$ can generically be 
avoided by the hypersurface equations $f=0$ and $g = 0$. 
For special values of the
coefficients $\{a_m \}$ and  $\{b_n \}$ in 
eq.\,\,(\ref{CYeq}) the corresponding Calabi--Yau varieties 
$\widehat{X}_f$ and 
 $\widehat{X}_g$ 
may of course be singular.

{}From now on  we want to restrict the types of singularities 
of $\overline{X}_f \subset \IP_\D$ under consideration and  demand  
that all $2$-dimensional faces $\th^\circ$  of the dual polytope
$\D^\circ$ are either
unimodular triangles (i.e. spanned by a subset of a lattice basis) or 
parallelograms of minimal volume, whose two triangulations hence are 
unimodular. This implies that all nonisolated singularities of 
$\IP_\D$ are $1$-parameter families of 
toric conifold singularities defined by an equation 
$u_1u_2 - u_3 u_4 =0$. These families $\IT_\th$ of singularties 
are in one-to-one 
correspondence to $1$-dimensional faces  $\th \subset \D$ such that the dual 
face $\th^\circ \subset \D^\circ$ is a parallelogram.

The morphism $\IP_{\S^\circ} \to \IP_\D$ induces a small 
crepant resolution $\widehat{X}_f \to \overline{X}_f$ which replaces 
every conifold point in  $\overline{X}_f$ by a copy of  $\IP^1$. 
We remark that every $1$-parameter family  $\IT_\th \subset \IP_\D$ 
(the $2$-dimensional dual face   $\th^\circ \subset \D^\circ$ 
is a parallelogram) 
of toric conifold singularities has exactly $l(\th) -1$ distinct 
intersection points with the generic hypersurface 
$\overline{X}_f \subset \IP_\D$. In order to 
analyse the possibility of smoothing $\overline{X}_f $ by a flat deformation 
we apply the following criterion of Namikawa 
(we formulate it in a simplified version): 
\\[7pt]
{\bf Theorem} (see \cite{Smoothing}): Let $X$ be a Calabi-Yau $3$-fold 
with $n$ isolated conifold singularities $$\{p_1, \ldots, p_k \} 
= Sing\, X$$ and $f\; : \; Z \to X$ be a small resolution 
of these singularities such that $C_i := f^{-1}(p_i) \cong 
\IP^1$ and $f$ an isomorphism over 
$X \setminus \{ p_1, \ldots, p_k \}$. Then $X$ can be deformed to a smooth 
Calabi-Yau $3$-fold if and only if the homology classes 
$[ C_i] \in H_2(Z, \IC)$ satisfy  a linear relation 
\[ \sum_{i =1}^k \alpha_i [ C_i ] = 0, \]
where $\alpha_i \neq 0$ for all $i$. It is easy to see 
that the last condition is equivalent to the fact that 
the subspace in  $H_2(Z, \IC)$ generated by the homology classes 
 $[ C_1 ], \ldots, [C_k ]$ coincides with the subspace 
generated by $\{ [ C_1 ], \ldots, [C_k ]\} \setminus \{ [C_i]\}$ 
for all $i =1, \ldots, k$.

In our situation we can choose $Z$ to be $\widehat{X}_f$. 
Let $P(\D)$ be the set of all $1$-dimensional faces of $\Delta$ 
such that the dual $2$-face $\th^\circ$  is a parallelogram of minimal 
volume. We set $k_\th = l(\theta) -1$. Then a generic  Calabi-Yau 
hypersurface $\overline{X}_f \subset \IP_\D$ contains
exactly
\[ k = \sum_{\th \in P(\D)} k_\th \]
conifold points. Let $\{ v_1, \ldots, v_l \}$ be the set of all vertices 
of the dual polytope $\D^\circ$. Then it can be shown that the 
homology group $H_2(\widehat{X}_f, \IQ)$ can be identified with 
the subgroup $R_{\D^\circ} \subset \IQ^l$ consisting of rational vectors 
$(\lambda_1, \ldots, \lambda_l) \in \IQ^l$ such that 
\[ \sum_{i =1}^l \lambda_i v_i =0. \]
In order to determine the homology class $[C_i] \in 
H_2(\widehat{X}_f, \IQ)$ we first remark that locally for each conifold 
point $p_i \in  \overline{X}_f $ there exist exactly two different 
small resolutions $f_i$ and $f_i'$ of $p_i$. The corresponding 
homology classes of exceptional curves differ by their signs 
$[f_i^{-1}(p_i)] = - [ (f_i')^{-1}(p_i)]$. The next step is to see that 
for any $\th \in P(\Delta)$ all $k_\th$ conifold points in the intersection 
$\IT_\th \cap  \overline{X}_f $ define (up to signs) the same homology 
class in  $R_{\D^\circ} \cong H_2(\widehat{X}_f, \IQ)$ coming from 
the linear relation 
\[ \rho_{\theta}\;: \; v_i +  v_j - v_s - v_r =0 ,\]
where $v_i, v_j, v_s, v_r$ are vertices of the parallelogram $\th^\circ$ 
and $[v_i,v_j ]$, $[v_s, v_r]$ are the two diagonals of  $\th^\circ$. 
Therefore, the homology classes
$[C_1], \ldots, [C_k] \in  H_2(\widehat{X}_f, \IQ)$ coincide (up to signs)
with the elements $\rho_\th$ for $\th\in P(\D)$. Each element $\pm\rho_\th$ 
appears $k_\th = l(\th) -1$ times in the sequence  $[C_1], \ldots, [C_k]$. 
Thus, the smoothing criterion of 
Namikawa can be formulated for   $ \overline{X}_f $ as follows: 

{\bf Smoothing criterion:} Under the above assumtion on $\D^{\circ}$, 
a generic Calabi-Yau hypersurface 
 $ \overline{X}_f  \subset \IP_\D$ is smoothable to a Calabi-Yau 
$3$-fold $Y$ by a flat deformation 
if and only if for any $1$-dimensional face $\th \in P$ such that 
$k_\th =1$ the 
element $\rho_\th$ is a linear combination of the remaining 
elements $\rho_{\th'}$ with $\th' \in P(\D), \;\; \th' \neq \th$.

Using the classification of $4$-dimensional reflexive 
polytopes, one can show that there exist exactly  $198849$ 
reflexive polytopes $\D$ such that all $2$-dimensional 
faces of the dual polytope $\D^\circ$ are either basic triangles, 
or parallelograms of minimal volume. Let $p:= |P(\D)|$ 
and $l$ be the number of vertices of $\D^\circ$. We define 
the matrix $\Lambda(\D)$ of size $p \times l$ whose 
rows are coefficients of the linear relation $\rho_\th$ 
$( \th \in P(\Delta))$. Then a 
generic Calabi-Yau hypersurface $ \overline{X}_f$ is smoothable by 
a flat deformation if and only if for all $\th \in P(\D)$ 
such that $k_\th =1$ 
the removing of the corresponding line $\L_\th(\D)$ 
from  $\L(\D)$ does not reduce the
rank of the matrix. 
This smoothing condition reduces the number of relevant polytopes from 198849
to 30241 as detailed according to Picard numbers $h^{1,1}$ in 
{\it table~\ref{TabPoly}}.

\begin{figure}[t]

\def\lab#1)#2#3{\put#1){\makebox(0,0)[#2]{\small #3}}}
\def\putlin#1,#2,#3,#4,#5){\put#1,#2){\line(#3,#4){#5}}} 
\def\putvec#1,#2,#3,#4,#5){\put#1,#2){\vector(#3,#4){#5}}}
	\def\Xlab#1 {\put(#1,#1){\drawline(-.5,.5)(.5,-.5)}}
	\def\Ylab#1 {\put(-#1,#1){\drawline(.5,.5)(-.5,-.5)}}

\input ccy.data
							\def\MaxCHC{16}
\vbox to 117pt	{\vss	\unitlength=2.7pt		\def\MaxNUM{130}
\begin{center}	\def\putlin#1,#2,#3,#4,#5){\put#1,#2){\line(#3,#4){#5}}}
		\def\putvec#1,#2,#3,#4,#5){\put#1,#2){\vector(#3,#4){#5}}}
	\def\Xlab#1 {\put(#1,#1){\drawline(-.5,.5)(.5,-.5)}}
	\def\Ylab#1 {\put(-#1,#1){\drawline(.5,.5)(-.5,-.5)}}

\unitlength=3.7pt		\def\CSH{.6} \def\CSC{.8}	
\BP(80,25)(25,2)	\putvec(0,0,1,0,130)	\putvec(0,0,0,1,\MaxCHC)
	\put(1,\MaxCHC){$h_{11}$}	\put(128,-3.2){$h_{21}$}
	\def\Xlab#1 {\put(#1,0){\drawline(0,-.5)(0,.5)}}
	\def\Ylab#1 {\put(0,#1){\drawline(.5,0)(-.5,0)}}
	\Xlab10 \Xlab20 \Xlab30 \Xlab40 \Xlab50 \Xlab60 \Xlab70 \Xlab80 
	\Xlab90	\Xlab100 \Xlab110 \Xlab120 \Xlab130 		\Ylab10

\def\h#1.#2.{ \ifnum \MaxCHC < #1 \else \put(#2,#1){\blue \circle{\CSC}} \fi }
\ConifoldCY

\def\h#1.#2.{ \ifnum \MaxNUM < #1 \else
	\ifnum \MaxCHC < #1 \else \put(#2,#1){\black\circle*{\CSH}} \fi
	\ifnum \MaxCHC < #2 \else \put(#1,#2){\black\circle*{\CSH}} \fi \fi}
\HypersurfCY
\EP\end{center}
					\stepcounter{figure}	\makeatletter
\immediate\write\@auxout{\string\newlabel{Conifolds}{\thefigure}}\makeatother
\vspace{8pt}
\centerline{Fig. \thefigure: Smoothed conifolds (circles) and toric hypersurfaces 
	(dots) with $h^{1,1}\le\MaxCHC$, $h^{2,1}\le130$.}
\vspace{1pt}}
\end{figure}

\begin{table}
\medskip
\centerline{\small	\begin{tabular}{||l|rrrrrrrrrrrrr||}\hline\hline
Picard number&1& 2  &3  &4  &5 &6&7&8  &9 &10 &11 &12 &15\\\hline\hline
Polytopes &8871&43080&74570&50863& 17090& 3540& 646& 124& 41& 17&2&4&1\\\hline
Smoothable& 210& 3470&11389&10264&3898&815& 140&35&9&8&1&1&1\\\hline\hline
			\end{tabular}}
\caption{Numbers of polytopes for conifold	\label{TabPoly}
	Calabi--Yau spaces with Picard number $h^{1,1}$.  }
\end{table}

The Hodge numbers of the smoothed  Calabi--Yau 3-fold $Y$  can be 
computed by the well-known formula (see e.g. \cite{DDP})
\BE	\label{Hodge}
	h^{1,1}(Y)=h^{1,1}(\widehat{X}_f) -\hbox{\tt rk},\qquad
	h^{1,2}=h^{1,2}(\widehat{X}_f) +\hbox{\tt dp}-\hbox{\tt rk}
\EE
where $\hbox{\tt rk}$ is the rank of the matrix $\L(\D)$ 
of linear relations and 
$\hbox{dp} = k = \sum_{\th \in P(\D)} k_\th$ 
denotes the number of double points in  $ \overline{X}_f$.
For the smoothable cases they are listed in {\it table~\ref{TabHodge}} and 
displayed as circles over the background of toric hypersurface data in
{\it fig.~\ref{Conifolds}}. The complete data, which was computed 
using the software package PALP \cite{PALP}, is available on the 
internet~\cite{web}.

\begin{table}	
\centerline{\small	\begin{tabular}{||l|r|l||}\hline\hline
$h^{1,1}$	& $\#(\D)$ & $h^{2,1}$	\\\hline\hline
1&210	&25,28--41,45,47,51,53,55,59,61,65,73,76,79,89,101,103,129\\\hline
2&3470	&26,28--60,62--68,70,72,74,76,77,78,80,82--84,86,88,90,96,100,102,112,116,128\\\hline
3&11389	&25,27--73,75--79,81,83,85,87,89,91,93,95,99,101,103,105,107,111,115\\\hline
4&10264	&24,28,30--76,78--82,84,86,88--98,100,102,104,106,112\\\hline
5&3898	&27,29,30--83,85--93,97\\\hline
6&815	&28,30--32,34--56,58--70,72--76,80,82\\\hline
7& 140	&27,29--31,33--35,37--41,43,45,47,49--51,53,55,57,59,61,62,64,76\\\hline
8&35	&30,32--34,36,38,40,42,44,52\\\hline
9&9	&31,33,37\\\hline
10&8	&26,30,34,36\\\hline
11&1	&27\\\hline
12&1	&28\\\hline
15&1	&23\\\hline\hline				\end{tabular}}
\caption{Hodge data $(h^{1,1},h^{2,1})$ for the 30241 smoothed conifold
	Calabi--Yau spaces.	\label{TabHodge}}
\end{table}

The enumeration of the polytopes $\D$ and of the Hodge data is, of course, 
only the first step and further work is required to compute 
the additional data like intersection form, Chern classes and quantum
cohomology of $Y$. This program will be initiated for the 210
examples of $Y$ with Picard number 1 in section 4.

\section{Mirror families} 

Let us discuss the explicit construction of mirrors of Calabi-Yau $3$-folds 
$Y$ from the previous section. As we have seen $Y$ is obtained 
from $\widehat{X}_f$ by a conifold transition. It turns out that the 
mirror family $Y^*$ of $Y$ is obtained also by a conifold transition 
from the mirrors $\widehat{X}_g = (\widehat{X}_f)^*$. 

For this we have to specialize the generic family of Laurent polynomials 
\[
	g = \sum_{n \in \D^\circ \cap N} b_n t^n 
\]
to a special one 
by imposing additional conditions on the coeffients
$b_n$: for every  $\th \in P(\Delta)$ 
 we demand 
\BE \label{special}
	 b_{v_i} b_{v_j} = b_{v_r} b_{v_s}, 
\EE
where  $v_i, v_j, v_r, v_s$ are vertices 
of the parallelogram $\th^\circ$ satisfying the equation 
\[ v_i + v_j = v_s + v_r. \]
We denote the specialized Laurent polynomial by $\tilde{g}$. 

Our main observation is that the hypersufaces 
$\widehat{X}_{\tilde{g}}$ from the specialized family 
should have the same number  $k = \sum_{\theta \in P(\Delta)} k_{\theta}$ 
of conifold singularities so that 
we could  consider  $\widehat{X}_{\tilde{g}}$ as a flat conifold degeneration 
of smooth Calabi-Yau $3$-folds  $\widehat{X}_{{g}}$.

Let us explain this observation in more detail. For any $\theta \in P(\Delta)$
we can choose a basis $e_1, \ldots, e_4$ of the lattice $M$ and the 
dual basis $e^{\circ}_1, \ldots, e^{\circ}_4$ of the dual lattice $N$ 
in such a way that $-e_3$ and $-e_3 - k_{\theta}e_4$ are vertices 
of the $1$-dimensional face $\theta \subset \Delta$ and 
$$e_3^{\circ}, \;\; e_3^{\circ} + e_1^{\circ},\;\;  
e_3^{\circ} + e_2^{\circ}, \;\;  e_3^{\circ} + e_2^{\circ} + e_1^{\circ}$$
are vertices of the dual (parallelogram) face 
$\theta^{\circ} \subset \Delta^{\circ}$. We put  $w_j:=-e_3 - j e_4 
\in  \theta \cap M$  $(0 \leq j \leq k_{\theta})$. 
Then 
\[ \theta \cap M = \{ w_0, \ldots, w_{ k_{\theta}} \} . \]
Every pair of lattice points $w_j, w_{j-1}$ $( 1 \leq j \leq k_{\theta})$ 
generates a $2$-dimensional cone $\sigma_j$ in the fan $\Sigma$. 
Since $e_1, e_2, w_j, w_{j-1}$ is a $\IZ$-basis of $M$ the cone $\sigma_j$ 
defines an affine open torus invariant subset $U_{\sigma_j} \subset 
\IP_{\Sigma}$ 
such that $U_{\sigma_j} \cong (\IC^*)^2 \times \IC^2 $. 
 Denote $w^{(j)}_3:=
(j-1)e_3^{\circ} - e_4^{\circ}$ and  $w^{(j)}_4 :=
-je_3^{\circ} + e_4^{\circ}$. Then $e_1^{\circ}, e_2^{\circ}, 
w^{(j)}_3, w^{(j)}_4$ is a $\IZ$-basis 
of $N$ dual to  $e_1, e_2, w_j, w_{j-1}$. 
We use this basis in order to define the local 
coordinates $t_1,t_2, t_3^{(j)},t_4^{(j)}$ 
on $U_{\sigma_j}$.  Then the  equation of 
$\widehat{X}_{\tilde{g}}$ in $U_{\sigma_j}$ can be written as follows
\[ \tilde{g}_j(t) = 
b  + b t_1 + bt_2 + b t_1 t_2 + 
t_3^{(j)} t_4^{(j)} \sum_{n \in A_j} b_n t^{n} \in 
\IC[t_1^{\pm 1},t_2^{\pm 1}, t_3^{(j)}, t_4^{(j)}], \]
where $A_j \subset N$ is a finite number of lattice points $n = 
(n_1,n_2,n_3, n_4) \in N$ 
satisfying the conditions $n_3 \geq 0$, $n_4 \geq 0$.   
Therefore, the polynomial 
\[ \tilde{g}_j(t) = 
b(1+t_1) (1+t_2)   + b_0  t_3^{(j)} t_4^{(j)} + 
t_3^{(j)} t_4^{(j)} \sum_{n \in A_j \setminus \{0\}} b_n t^{n} \]
has a conifold singularity at 
the point $q_j = (-1,-1,0,0) \in U_{\sigma_j}$ 
$( 1 \leq j \leq k_{\theta})$.  By repeating this computation for
every $1$-dimensional face $\theta \in P(\Delta)$, we obtain  
 $k = \sum_{\theta \in P(\Delta)} k_{\theta}$ conifold points
in  $\widehat{X}_{\tilde{g}}$.

{\bf Remark.} Unfortunately, these above arguments do not show 
that we have found {\bf all} singular points of  $\widehat{X}_{\tilde{g}}$. 
We hope that this is true in many cases.

Now we are going to obtain the mirrors $Y^*$ using small resolutions of 
singularities of  $\widehat{X}_{\tilde{g}}$. 
By a result of Smith, Thomas and Yau \cite{STY} (Theorem 2.9), in order that
$\widehat{X}_{\tilde{g}}$ admits a projective small resolution, the homology 
classes $[L_1], \ldots,[L_k] \in H_3 (\widehat{X}_{{g}}, \IZ)$ 
of the vanishing cycles $L_i \cong S^3$ must satisfy a linear relation 
\[ \sum_{i =1}^k c_i [L_i] =0, \; \; c_i \neq 0 \;\; \forall i. \]
This condition can be considered as ``mirror'' to the criterion of Namikawa.

In the case when $h^{1,1}(Y) =1$ the specialization equations 
(\ref{special}) show that we can put $b_0=1$ and  $b_n =-z$ $\forall n 
\neq 0)$, so that we obtain a one-parameter family 
of Laurent polynomials $\tilde{g}$ 
depending only on $z$.


\begin{table}[h]		
\newcounter{ncy}   \setcounter{ncy} 0   \def\NT{\stepcounter{ncy}{\#\thency}}
\def\ETL{\\\hline}
\def\bt{\hspace{-6pt}\begin{tabular}l}  \def\et{\end{tabular}\hspace{-7pt}}
\vbox {\vspace*{-1pt}
\begin{center}{\footnotesize
		\hspace*{-9mm}
\begin{tabular}{||r||c|c|c||c||c|c	
				||}\hline\hline   
	& $h_{12}$& $H^3$& $c_2H$& $c_3$ & 
			$\!\!N_\D)\!\!$ & $\!\!N_{\PP_0}\!\!$
							\\\hline\hline
\NT&\bf25 &79 &94 &-48&    1&1\ETL	
\NT&\bf28 &99 &102 &-54&   1&1\ETL	
\NT&  28 &104 &104 &-54&  2&1\ETL	
\NT&\bf29 &74 &92 &-56& 2&2\ETL	
\NT&  29 &88 &100 &-56&   1&1\ETL	
\NT&  29 &93 &102 &-56&   4&3\ETL	
\NT&  29 &98 &104 &-56&   1&1\ETL	
\NT&\bf30 &98 &104 &-58&5&4\ETL 
\NT&  30 &103 &106 &-58& 2&2\ETL	
\NT&  30 &108 &108 &-58& 4&2\ETL	
\NT&\bf31 &78 &96 &-60&    1&1 \ETL	
\NT&  31 &83 &98 &-60&    2&2 \ETL	
\NT&  31 &98 &104 &-60&   6&4 \ETL	
\NT&  31 &103 &106 &-60&  3&1 \ETL	
\NT&  31 &108 &108 &-60&  1&1 \ETL	
\NT&  31 &118 &112 &-60&  5&2 \ETL	
\NT&  31 &124 &112 &-60&  2&1 \ETL	
\NT&\bf32 &83 &98 &-62&    2&2 \ETL	
\NT&  32 &98 &104 &-62&   3&1 \ETL	
\NT&  32 &108 &108 &-62&  4&3 \ETL	
\NT&  32 &113 &110 &-62&  1&1 \ETL	
\NT&  32 &118 &112 &-62&  4&3 \ETL	
\NT&\bf33 &78 &96 &-64&    1&1 \ETL	
\NT&  33 &97 &106 &-64&   1&1 \ETL	
\NT&  33 &108 &108 &-64&  4&1 \ETL	
\NT&\bf34 &97 &106&-66& 3&3 \ETL	
\NT&  34 &102 &108 &-66&  6&3 \ETL	
\NT&  34 &123 &114 &-66&  1&1 \ETL	
\NT&\bf35 &87 &102 &-68& 1&1 \ETL	
\NT&  35 &92 &104 &-68&   7&5 \ETL	
\NT&  35 &97 &106 &-68&   5&3 \ETL	
\NT&  35 &102 &108 &-68&  8&4 \ETL	
\NT&  35 &112 &112 &-68&  13&3 \ETL	
\NT&\bf36 &82 &100 &-70&  1&1 \ETL	
\end{tabular}
~~
~~
\begin{tabular}{||r||c|c|c||c||c|c||c||}\hline\hline   
	\#& $h_{12}$& $H^3$& $c_2H$& $c_3$ 
		& $\!\!N_\D\!\!$ & $\!\!N_{\PP_0}\!\!$&${PF}$\\\hline\hline
\NT&  36 &92 &104 &-70&  5&3 &\ETL	
\NT&  36 &107 &110 &-70&  16&5 &\ETL	
\NT&\bf37 &117 &114 &-72&  12&1 &(\ref{X37})\ETL	
\NT&\bf38 &102 &108 &-74&  2&1 &\ETL	
\NT&\bf39 &96 &108 &-76&  2&1 &\ETL	
\NT&  39 &152 &116 &-76&  1&1 &(\ref{X39})\ETL	
\NT&\bf40 &91 &106 &-78&  2&1 &\ETL	
\NT&\bf41 &86 &104 &-80&  2&1 &\ETL	
\NT&  41 &116 &116 &-80&  13&1 &(\ref{X41})\ETL	
\NT&\bf45 &144 &120 &-88&  2&2&$\!\!$(\ref{new-4})$\,,\,$\bf214$\!\!$\ETL
\NT&\bf47 &144 &120 &-92&  2&1& \bf289\ETL	
\NT&  47 &176 &128 &-92&  3&1&(\ref{new-7})\ETL	
\NT&\bf51 &168 &132 &-100& 1&1&\bf218\ETL	
\NT&  51 &200 &140 &-100& 3&2&$\!\!$(\ref{new-9})$\,,\,$(\ref{new-10})$\!\!$
					\ETL	
\NT&\bf53 &168 &132 &-104& 2&1& \bf287 \ETL	
\NT&  53 &232 &148 &-104& 4&1&(\ref{new-2})\ETL	
\NT&\bf55 &136 &124 &-108& 1&1&\bf209\ETL	
\NT&\bf59 &24 &72 &-116&2&1   &\bf29\ETL	
\NT&  59 &28 &76 &-116& 3&1   &\bf26\ETL	
\NT&  59 &32 &80 &-116& 4&1   &\bf42\ETL	
\NT&\bf61 &20 &68 &-120& 1&1  &\bf25\ETL	
\NT&  61 &36 &84 &-120&  4&1  &\bf185\ETL	
\NT&\bf65 &16 &64 &-128& 1&1  &\bf3\ETL	
\NT&  65 &44 &92 &-128&  1&1 &(\ref{new-?})\ETL
\NT&\bf73& 9 &54 &-144& 1&1   &\bf4\ETL	
\NT&  73 &12 &60 &-144&  2&1  &\bf5\ETL	
\NT&  73 &32 &80 &-144&  1&1  &\bf10\ETL	
\NT&\bf76 &15 &66 &-150&  2&1 &\bf24\ETL	
\NT&\bf79 &48 &96 &-156&   1&1 &\bf11\ETL	
\NT&  79 &432 &192 &-156&   1&1 &\bf12\ETL	
\NT&\bf89 &8 &56 &-176&     2&1 &\bf6\ETL	
\NT&\bf101 &80 &128 &-200&  1&1 &\bf51\ETL	
\NT&\bf103 &648 &252 &-204& 1&1 &\bf8\ETL	
\NT&\bf129 &108 &156 &-256& 1&1 &\bf14\ETL	
\end{tabular}\vspace{-15pt}
}\hspace*{-12mm}
\end{center}
\caption{
	Topological data with multiplicities $N_\D$ of polytopes
	and  $N_{\PP_0}$ of principal periods.
\hspace*{47pt}	The last column, denoted $PF$, refers to the Picard--Fuchs 
	operator, either by equation 
\hspace*{47pt}	number (in parentheses) or (in boldface) by the reference
	number in the tables of \cite{DucoList}.
}	\label{TabPic1}
}
\vspace{-9pt}
\end{table}

\section{One-parameter manifolds}
\vspace{-9pt}
We now focus on the list of the 210 polytopes that lead to one-parameter
families with smoothable conifold singularities. According to a theorem 
by Wall \cite{Wall} the diffeomorphism type of a Calabi--Yau is completely 
characterized by its Hodge numbers, intersection ring and second Chern class.
For Picard number one the latter two amount to the tripple intersection
number $H^3$ and the number $Hc_2$. The resulting 68 different 
topological types are listed in     {\it table~\ref{TabPic1}}.

The entries in {\it table~\ref{TabPic1}} have been computed combinatorially 
with the formulas
\BE							\label{intersect}
	H^3={\rm Vol(\D)}/({\rm Ind})^3, \qquad 
	c_2\cdot H=(12\,|\partial\D\cap M|\,)	/{\rm Ind}
\EE
where $\rm Vol$ denotes the lattice volume, $|\partial\D\cap M|$ is the number
of boundary lattice points of $\D$ and $\rm Ind$ is the index of the
affine sublattice of $M$ that is generated by the vertices of $\D$. There
is one exception to this rule, namely the convex hull of the Newton
polytope of%
\footnote{~This is one of the five polytopes for which a facet of $\D$ has
		an interior point.}
\BE	
	g=t_1 + t_1t_2^2 + t_1t_2^{2}t_3^{4} + t_1t_2^{2}t_4^{4} 
	+ t_1^{-3}t_2^{-2}t_3^{-4} + t_1^{-3}t_2^{-2}t_4^{-4} 
	+ t_1^{-3}t_2^{-2},
\EE
which is one of the two polytopes that lead to the entry \#65 with $h_{12}=89$
in the table. For this variety the divisor $H$ has multiplicity two so 
that effectively $\rm Ind$ has to be doubled in eq.~(\ref{intersect}).

A glance at {\it fig.\,\ref{Conifolds}} shows that we constructed a 
surprisingly rich new set as compared to toric hypersurfaces, and also
when compared to other know constructions of one-parameter examples
\cite{kkrs,Grass,Flag,MorePic1}.
Observe that in all instances in {\it table~\ref{TabPic1}} with given 
Hodge numbers each of the intersection 
numbers $H^3$ and $c_2\!\,\cdot\!H$ uniquely determines the other.

For the computation of the quantum cohomology we start with the computation
of the principal period	
\BE			\label{Pint}
     \PP_0(z)=\oint \frac{dt_1}{t_1}\frac{dt_2}{t_2}\frac{dt_3}{t_3}
	\frac{dt_4}{t_4} \Biggl(1-z\sum_{m\in Vert(\D^\circ)} 
t^m\Biggr)^{-1}\!.
\EE
where for $h_{11}=1$ all coefficients of the relevant non-constant
monomials can be set to $b_j=-z$ because mirror symmetry amounts to the
restriction $b_{j_1}b_{j_2}=b_{j_3}b_{j_4}$ if the corresponding vertices
form a parallelogram with $v_{j_1}+v_{j_2}=v_{j_3}+v_{j_4}$ and 
$h^{2,1}(\widehat{X}_g)-rk=1$ implies that rescalings of the homogeneous 
coordinates can be used to identify the corresponding coefficients with 
the complex structure modulus $-z$. The function 
$\PP_0(z)$ is the unique regular
power series solution in the kernel of the Picard--Fuchs operator
\BE	\label{PFO}
	{\cal O}=\th^4+\sum_{n=1}^d ~z^n~\sum_{i=0}^4 c_{ni}\th^i,~\qquad
	\th=z\frac d{dz}.
\EE
This operator can then be used to compute the periods with logarithmic
singularities and the instanton numbers via the mirror map as explained, 
for example, in \cite{CoKa}.
Our method for the computation of $\cal O$ is direct evaluation of the 
period by expansion of the last term in eq. (\ref{Pint})
in $z$ up to (at least) degree $5d$ and determination
of the coefficients $c_{ni}$ from ${\cal O}\,(\PP_0)=0$
for the ansatz eq. (\ref{PFO}).

We have computed all Picard--Fuchs operators for the manifolds with
$h_{12}\ge45$, which are the cases \#44 -- \#68 in {\it table~\ref{TabPic1}}. 
They have been determined independently by Duco van Straten and
Gert Almkvist \cite{Duco}.
Most of these operators were known before \cite{DucoList}. Here we 
only list three examples that are needed, 
in addition to eqs. (\ref{44a}), (\ref{48a}) and (\ref{48b}) below,
as references in {\it table~\ref{TabPic1}}:	

\def\cy #1 #2 #3 {$\mathbf{X^{#1}_{#2,#3}}$}	\let\nn=\nonumber
\def\top h12=#1 E=#2 H^3=#3 c2H=#4 {\normalsize~\\[-15pt]
	\rm\#\NP: \cy #1 #3 #4 \footnotesize}
\def\NLQ{\\\nn&&\hspace{18pt}}	\let\NLF=\NLQ

\medskip
\def\NP{58}	\top h12=65 E=-128 H^3=44 c2H=92
\BEA&&	\textstyle						\label{new-?}
\th^4 - 4~z~\th (\th+1) (2 \th+1)^2 
 - 32~z^2~(2 \th+1) (2 \th+3) (11 \th^2+22 \th+12) 
\NLQ - 1200~z^3~(2 \th+1) (2 \th+3)^2 (2 \th+5) 
 - 4864~z^4~(2 \th+1) (2 \th+3) (2 \th+5) (2 \th+7) 
\EEA

\def\NP{50}	\top	h12=53 E=-104 H^3=232 c2H=148
\BEA&&	\textstyle						\label{new-2}
\th^4 - \frac{2}{29}~z~(1318 \th^4+2336 \th^3+1806 \th^2+638 \th+87) 
\NLQ - \frac{4}{841}~z^2~(90996 \th^4+744384 \th^3+1267526 \th^2+791584 \th+168345) 
\NLQ + \frac{100}{841}~z^3~(34172 \th^4+77256 \th^3-46701 \th^2-110403 \th-36540) 
\NLF + \frac{10000}{841}~z^4~(2 \th+1) (68 \th^3+1842 \th^2+2899 \th+1215) 
 - \frac{5000000}{841}~z^5~(\th+1)^2 (2 \th+1) (2 \th+3) 		
\EEA

\def\NP{46}	\top	h12=47 E=-92 H^3=176 c2H=128
\BEA&&	\textstyle						\label{new-7}
\th^4 - \frac{4}{11}~z~(432 \th^4+624 \th^3+477 \th^2+165 \th+22) 
\NLQ + \frac{32}{121}~z^2~(12944 \th^4+4736 \th^3-15491 \th^2-12914 \th-2860) 
\NLQ - \frac{80}{121}~z^3~(10688 \th^4-114048 \th^3-159132 \th^2-83028 \th-15455) 
\NLF - \frac{51200}{121}~z^4~(2 \th+1) (4 \th+3) (76 \th^2+189 \th+125) 
 + \frac{2048000}{121}~z^5~(2 \th+1) (2 \th+3) (4 \th+3) (4 \th+5) 
\EEA
\normalsize

The calculations become quite expensive when the number of vertices of
$\D^\circ$ becomes large, as is mostly the case for manifolds with small
$h_{12}$. Nevertheless we could determine, so far, the operators for three 
more examples:

\noindent
\\[3pt]\#43: \cy 41 116 116 : \footnotesize
\BEA	&&\!\!\!\!\!\!\!\!\!
	\th^4 + \frac{2}{29}~z~\th (24 \th^3-198 \th^2-128 \th-29) 
- \frac{4}{841}~z^2~(44284 \th^4+172954 \th^3+248589 \th^2+172057 \th+47096) 
\nn\\&& - \frac{4}{841}~z^3~(525708 \th^4+2414772 \th^3
		+4447643 \th^2+3839049 \th+1275594) 
\nn\\&& - \frac{8}{841}~z^4~(1415624 \th^4+7911004 \th^3
		+17395449 \th^2+17396359 \th+6496262) 
\nn\\&& - \frac{16}{841}~z^5~(\th+1) (2152040 \th^3+12186636 \th^2
		+24179373 \th+16560506) 
\nn\\&& - \frac{32}{841}~z^6~(\th+1)(\th+2)(1912256\th^2+9108540\th+11349571) 
\nn\\&& - \frac{10496}{841}~z^7~(\th+1) (\th+2) (\th+3) (5671 \th+16301) 
	- \frac{24529152}{841}~z^8~(\th+1) (\th+2) (\th+3) (\th+4) \label{X41}
\EEA
\normalsize
\\[3pt]\#40: \cy 39 152 116 : \footnotesize
\BEA	&&\!\!\!\!\!\!\!\!\!
	\th^4 - \frac{1}{19}~z~(4333 \th^4+6212 \th^3+4778 \th^2+1672 \th+228)
\nn\\&&	+ \frac{1}{361}~z^2~(4307495 \th^4+7600484 \th^3
		+6216406 \th^2+2802424 \th+530556)
\nn\\&& - \frac{1}{361}~z^3~(93729369 \th^4+213316800 \th^3
		+236037196 \th^2+125748612 \th+25260804) 
\nn\\&&	+ \frac{4}{361}~z^4~(240813800 \th^4+778529200 \th^3
		+1041447759 \th^2+631802809 \th+138510993)
\nn\\&&	- \frac{1636}{361}~z^5~(\th+1) (2851324 \th^3
		+10035516 \th^2+11221241 \th+3481470)
\nn\\&&	+ 6022116~z^6~(\th+1) (\th+2) (2 \th+1) (2 \th+5) 	\label{X39}
\EEA
\normalsize
\\[3pt]\#37: \cy 37 117 114 :\footnotesize	
\BEA	&&\!\!\!\!\!\!\!\!\!
	\th^4 - \frac{1}{13}~z~\th (56\th^3+178 \th^2+115 \th+26) 
	- \frac1{169}~z^2~(28466\th^4+109442\th^3+165603\th^2+117338\th+32448)
\nn\\&&	- \frac{1}{169}~z^3~(233114 \th^4+1257906 \th^3
		+2622815 \th^2+2467842 \th+872352)
\nn\\&&	- \frac{1}{169}~z^4~(989585 \th^4+6852298 \th^3
		+17737939 \th^2+19969754 \th+8108448) 
\nn\\&&	- \frac1{169}~z^5~(\th+1) (2458967\th^3+18007287\th^2+44047582\th
	+35386584)
\nn\\&&	- \frac{9}{169}~z^6~(\th+1) (\th+2) (393163 \th^2+2539029 \th
	+4164444) 
\nn\\&&	- \frac{297}{169}~z^7~(\th+1) (\th+2) (\th+3) (8683 \th+34604) 
	- \frac{55539}{13}~z^8~(\th+1) (\th+2) (\th+3) (\th+4) \label{X37}
\EEA
\normalsize

The largest degree of a coefficient that we computed so far is degree 65, 
which we did for the conifolds \#17 and \#28, so that in these
cases the Picard--Fuchs operators would have at least degree 14. Further 
results will be put on our data supplement web page at \cite{web} as they
become available.

\subsection{Fractional transformations and instanton numbers}

Even though we do not yet know the Picard--Fuchs operators in many cases
it can be seen already from the first terms in the power series expansion
of the principal period which polytopes will yield identical operators.
In addition to the degeneracy of up to 13 different polytopes yielding 
the same Picard--Fuchs operator we thus observe that the same diffeomorphism
type can yield up to 5 different Picard--Fuchs operators, as indicated
in {\it table~\ref{TabPic1}}. Among the operators that we know this
phenomenon occurs twice:

\normalsize	
~		\\[-7pt]For \#44: \cy 45 144 120 \ we find
\footnotesize
\BEA
	\th^4 \!\!\!\!&-&\!\! 2~z~(102 \th^4+204 \th^3+155 \th^2+53 \th+7) 
\nn			~+~ 4~z^2~(\th+1)^2 (396 \th^2+792 \th+311) \\
\!\!\!&-&\!\!784~z^3~(\th+1)(\th+2)(2\th+1)(2\th+5) \label{new-4}\label{44a}
\EEA
\normalsize	for 
$f_\D=\frac{t_1t_4}{t_3}+\frac{t_2t_4}{t_1}+\frac{t_1t_4}{t_2t_3}+t_1t_4
	+\frac{t_2}{t_1}+\frac{1}{t_1t_4}
        +\frac{1}{t_1}+\frac{t_1}{t_2t_3}+\frac{t_1}{t_2t_4}
	+\frac{t_3}{t_1t_4}+\frac{t_2t_3}{t_1t_4}
        +\frac{t_1}{t_2}+\frac{t_2t_3}{t_1}+t_1    
$
and	\footnotesize
\BEA
\!\!\!\!\th^4 \!\!\!\!&-&\!\! 2~z~(90 \th^4+188 \th^3+141 \th^2+47 \th+6) 
 ~-~ 4~z^2~(564 \th^4+1520 \th^3+1705 \th^2+934 \th+192) 
\nn\\ 	&-&\!\! 16~z^3~(2 \th+1) (286 \th^3+813 \th^2+851 \th+294) \label{44b}
 ~-~ 192~z^4~(2 \th+1) (2 \th+3) (4 \th+3) (4 \th+5)            
\EEA
\normalsize	for 
$f_\D=\frac{1}{t_1}+\frac{t_4}{t_1}+\frac{t_2t_4}{t_1}+\frac{t_2}{t_1}
	+\frac{t_2t_3}{t_1}+\frac{t_2t_3t_4}{t_1}
        +\frac{t_3t_4}{t_1}+\frac{t_3}{t_1}+\frac{t_1}{t_2t_4}
	+\frac{t_1}{t_2}+\frac{t_1}{t_4}
        +\frac{t_1}{t_3t_4}+\frac{t_1}{t_3}+\frac{t_1}{t_2t_3} .
$
\normalsize	\\[9pt]For \#48: \cy 51 200 140 \ the two
operators, with three respective polytopes, are
\footnotesize
\BEA
\th^4 \!\!\!\!&-&\!\! ~z~(113 \th^4+226 \th^3+173 \th^2+60 \th+8) 
\nn ~-~ 8~z^2~(\th+1)^2 (119 \th^2+238 \th+92) \\		 \label{new-9}
 \!\!\!\!&-&\!\! 484~z^3~(\th+1) (\th+2) (2 \th+1) (2 \th+5)	 \label{48a} 
\EEA\normalsize
for $f_\D=t_1t_2+\frac{t_3}{t_2}+\frac{t_3t_4}{t_2}+\frac{t_2}{t_4}
	+\frac{1}{t_2}+\frac{t_2}{t_3}+\frac{t_4}{t_2}
     +\frac{t_3}{t_1t_2}+\frac{1}{t_1t_2}+\frac{t_4}{t_1t_2}
	+\frac{t_2}{t_3t_4}+\frac{t_2}{t_1t_3t_4} $
as well as for the Newton polytope of
$f_\D=\frac{t_2t_4}{t_1}+\frac{t_2}{t_1}+\frac{t_3}{t_1}+\frac{t_3t_4}{t_1}
	+\frac{1}{t_1}+\frac{t_4}{t_1}
        +\frac{t_2t_3}{t_1}+\frac{t_1}{t_4}+\frac{t_1}{t_2t_3}
	+\frac{t_1}{t_2t_3t_4}+\frac{t_1}{t_2t_4}+\frac{t_1}{t_3}
$
and \footnotesize
\BEA
\th^4 \!\!\!\!&-&\!\! z~(137 \th^4+258 \th^3+201 \th^2+72 \th+10) 
 ~+~ 4~z^2~(387 \th^4+1016 \th^3+1151 \th^2+642 \th+135) 		\nn
\\ \!\!\!\!&-&\!\! 4~z^3~(2 \th+1) (758 \th^3+2137 \th^2+2269 \th+820) 
 ~+~ 2000~z^4~(\th+1)^2 (2 \th+1) (2 \th+3) 	\label{48b} \label{new-10}
\EEA	\normalsize
for $f_\D=\frac{t_1}{t_2t_3}+\frac{t_1t_4}{t_2t_3}+\frac{t_2t_3}{t_1}
+\frac{t_2t_3}{t_1t_4}+\frac{t_2t_4}{t_1}
        +\frac{t_2}{t_1}+\frac{t_3}{t_1t_4}+\frac{1}{t_1}+\frac{1}{t_1t_4}
+\frac{t_1}{t_2}
        +\frac{t_1}{t_2t_4}+\frac{t_1t_4}{t_3} 
$.
\normalsize
\medskip

It is, of course, an interesting question whether the symplectic
Gromov--Witten invariants can give a finer classification than the 
diffeomorphism type. We do, however, not know a single example of such a
situation. We hence expect that the Picard-Fuchs operators 
(\ref{44a}-\ref{44b}) and (\ref{48a}-\ref{48b}) are related by rational
changes of variables that do not change the instanton numbers (cf. 
appendix A of \cite{kkrs}).
This is indeed the case.

For the diffeomorphism type \cy 45 144 120 \ the differential operator
(\ref{44a}) is transformed into (\ref{44b}) by the change of variables
\BE
	z\to\frac z{1+4z}
\EE
and the instanton numbers are, for both cases,
{\small
\BE
	n^{(0)}=\{3744,50112,1656320,77726016,4505800320,298578230016,
	21713403010176,
	\ldots\}.
\EE}%
For the diffeomorphism type \cy 51 200 140 \ the differential operator
(\ref{48a}) is transformed into (\ref{48b}) by the change of variables
\BE
	z\to\frac z{1-4z}
\EE
and the instanton numbers are, again for both cases,
{\small
\BE
	n^{(0)}=\{2600,25600,530000,15880000,584279000,24562482400,
	1132828485400,
	\ldots\}.
\EE}%
It will be interesting to check whether this phenomenon continues to hold
for the cases with smaller $h_{12}$ for which there are up to five 
different Picard--Fuchs operators for the same diffeomorphism type.


\RTransNotes
$\th^4 - 2 z (102 \th^4+204 \th^3+155 \th^2+53 \th+7)
      +  4 z^2 (\th+1)^2 (396 \th^2+792 \th+311)
      -   784 z^3 (\th+1) (\th+2) (2 \th+1) (2 \th+5)	$

$Y_{1,1,1}=-144\,{\frac {1}{ \left( 196\,z_{{1}}-1 \right)
	\left(4\,z_{{1}}-1 \right) ^{2}}}$

$q_{1} = z_{{1}}+50\,{z_{{1}}}^{2}+4205\,{z_{{1}}}^{3}+450414\,{z_{{1}}}^{4}+
55217174\,{z_{{1}}}^{5}+7376921468\,{z_{{1}}}^{6}+1046022991105\,{z_{{
1}}}^{7}+154944402713406\,{z_{{1}}}^{8}+23730680382255912\,{z_{{1}}}^{
9}+3731500623723105136\,{z_{{1}}}^{10}+599388854949993087752\,{z_{{1}}
}^{11}+97987073404363282282736\,{z_{{1}}}^{12}+
16256808045164134926962240\,{z_{{1}}}^{13}+
2731196146638410345114186240\,{z_{{1}}}^{14}+O \left( {z_{{1}}}^{15}\right)$

$\th^4 - 2 z (90 \th^4+188 \th^3+141 \th^2+47 \th+6)
  -  4 z^2 (564 \th^4+1520 \th^3+1705 \th^2+934 \th+192)
  -  16 z^3 (2 \th+1) (286 \th^3+813 \th^2+851 \th+294)
  -  192 z^4 (2 \th+1) (2 \th+3) (4 \th+3) (4 \th+5)	$  

$	Y_{1,1,1} = -144\,{\frac {1}{ \left( 4\,z_{{1}}+1 \right)  \left(
	192\,z_{{1}}-1 \right) }}			$

$q_{1} = z_{{1}}+46\,{z_{{1}}}^{2}+3821\,{z_{{1}}}^{3}+402290\,{z_{{1}}}^{4}+
48401686\,{z_{{1}}}^{5}+6342016004\,{z_{{1}}}^{6}+881668311809\,{z_{{1
}}}^{7}+128014664808290\,{z_{{1}}}^{8}+19215603047893288\,{z_{{1}}}^{9
}+2961051664483108944\,{z_{{1}}}^{10}+466079324081950468104\,{z_{{1}}}
^{11}+74659696858818473903760\,{z_{{1}}}^{12}+
12136728979396848918751808\,{z_{{1}}}^{13}+
1997811582059187625557496576\,{z_{{1}}}^{14}+O \left( {z_{{1}}}^{15} \right)$

Die Transformation vom ersten zum zweiten Operator geht wie folgt:
$	z_{1} \to \frac{ z_{1} }{ \left(1 + 4\,z_{1} \right) }	\qquad
	Y_{1,1,1} \to \frac{ Y_{1,1,1} }{ \left(1 + 4\,z_{1} \right)^4 }$
\\(Die Transformation von $Y_{1,1,1}$ entspricht einer Wahl des Schnittes 
von ${\cal L}^2$, cf. die Bemerkung nach (A.18) in hep-th/0410018).

$C_{1,1,1} = 144+3744\,q_{{1}}+404640\,{q_{{1}}}^{2}+44724384\,{q_{{1}}}^{3}+
4974869664\,{q_{{1}}}^{4}+563225043744\,{q_{{1}}}^{5}+64492942808736\,
{q_{{1}}}^{6}+7447697232494112\,{q_{{1}}}^{7}+865573069466565792\,{q_{
{1}}}^{8}+101098189285482631392\,{q_{{1}}}^{9}+11855085855044349284640
\,{q_{{1}}}^{10}+1394661536823040365810528\,{q_{{1}}}^{11}+
164511477390607751420765856\,{q_{{1}}}^{12}+
19449243716479906661441826912\,{q_{{1}}}^{13}+
2465158222457349901086047832576\,{q_{{1}}}^{14}+
312884551220283496908679444531776\,{q_{{1}}}^{15}+O\left({q_{1}}^{16}\right)$

$n^{(0)}=[3744,50112,1656320,77726016,4505800320,298578230016,21713403010176,
1690572391585344,138680643738597952,11855085854481123840,
1047829854863291033664,95203401267675493925632,
8852637103541150050724544]$

$\th^4 -  z (113 \th^4+226 \th^3+173 \th^2+60 \th+8)
  - 8 z^2 (\th+1)^2 (119 \th^2+238 \th+92)
  - 484 z^3 (\th+1) (\th+2) (2 \th+1) (2 \th+5)$  

$Y_{1,1,1} = -200\,{\frac {1}{ \left( 121\,z_{{1}}-1 \right)  \left( 
		4\,z_{{1}}+1 \right) ^{2}}}	$

$q_{1} = z_{{1}}+28\,{z_{{1}}}^{2}+1490\,{z_{{1}}}^{3}+99596\,{z_{{1}}}^{4}+
7592047\,{z_{{1}}}^{5}+629392604\,{z_{{1}}}^{6}+55310469032\,{z_{{1}}}
^{7}+5073426423380\,{z_{{1}}}^{8}+480885501488743\,{z_{{1}}}^{9}+
46777318238314276\,{z_{{1}}}^{10}+4646652664234583956\,{z_{{1}}}^{11}+
469646325905786507108\,{z_{{1}}}^{12}+48163878799360387851515\,{z_{{1}
}}^{13}+5000958294577194345017844\,{z_{{1}}}^{14}+O \left( {z_{{1}}}^{
	15} \right)$
\EndRTransNotes

\RTransNotes
$\th^4 - z (137 \th^4+258 \th^3+201 \th^2+72 \th+10)
  +  4 z^2 (387 \th^4+1016 \th^3+1151 \th^2+642 \th+135)
  -  4 z^3 (2 \th+1) (758 \th^3+2137 \th^2+2269 \th+820)
  +  2000 z^4 (\th+1)^2 (2 \th+1) (2 \th+3)$  

$Y_{1,1,1} = 200\,{\frac {1}{ \left( 125\,z_{{1}}-1 \right)  
	\left( 4\,z_{{1}}-1 \right) }}$

$q_{1} = z_{{1}}+32\,{z_{{1}}}^{2}+1730\,{z_{{1}}}^{3}+118884\,{z_{{1}}}^{4}+
9336047\,{z_{{1}}}^{5}+798159368\,{z_{{1}}}^{6}+72371363416\,{z_{{1}}}
^{7}+6851526896820\,{z_{{1}}}^{8}+670411915591911\,{z_{{1}}}^{9}+
67330155400411808\,{z_{{1}}}^{10}+6906089504745934564\,{z_{{1}}}^{11}+
720797237040176990548\,{z_{{1}}}^{12}+76337662450199976783675\,{z_{{1}
}}^{13}+8185884080167692792226416\,{z_{{1}}}^{14}+O \left( {z_{{1}}}^{
15} \right)$

Die Transformation vom ersten zum zweiten Operator geht wie folgt:\\
$z_{1} \to \frac{ z_{1} }{ \left(1 - 4\,z_{1} \right) }\qquad
	Y_{1,1,1} \to \frac{ Y_{1,1,1} }{ \left(1 - 4\,z_{1} \right)^4 }$

$C_{1,1,1} = 200+2600\,q_{{1}}+207400\,{q_{{1}}}^{2}+14312600\,{q_{{1}}}^{3}+
1016527400\,{q_{{1}}}^{4}+73034877600\,{q_{{1}}}^{5}+5305510715800\,{q
_{{1}}}^{6}+388560170494800\,{q_{{1}}}^{7}+28634333066754600\,{q_{{1}}
}^{8}+2120463844096346000\,{q_{{1}}}^{9}+157641571317395082400\,{q_{{1
}}}^{10}+11757031435553885965400\,{q_{{1}}}^{11}+
879180546255533783970200\,{q_{{1}}}^{12}+65891825571445476401052800\,{
q_{{1}}}^{13}+5349931173086603768195259600\,{q_{{1}}}^{14}+
429987891545197860293696313600\,{q_{{1}}}^{15}+O \left( {q_{{1}}}^{16}\right)$

$n^{(0)} = [2600,25600,530000,15880000,584279000,24562482400,1132828485400,
55926429785600,2908729552924600,157641571244360000,8833231732196758800
,508785038339251884800,29991727615587381156600]$
\EndRTransNotes

\RTransNotes
Ich habe ueber Deine Frage nach dem Transformationsverhalten der 
Yukawakopplungen nachgedacht. Ich war da ein wenig zu schnell mit dem 
Schreiben der E-mail letzte Woche. Schauen wir uns das erste Beispiel an:

$Y_{1,1,1}(z_{{1}})=-144\,{\frac {1}{ \left( 196\,z_{{1}}-1 \right)  \left(
	4\,z_{{1}}-1 \right) ^{2}}}$

Die Variablentransformation ist:
$z_{{1}}=\frac {x_{{1}}}{1+4\,x_{{1}}}$

Wenn man das in die Transformationsformel nach (A.18) einsetzt, findet man\\
$Y_{1,1,1}(x_{{1}}) = \frac{1_x}{l_z}{\frac {-144}{ \left( 4\,x_{{1}}+1
\right)^3  \left(192\,x_{{1}}-1 \right) }}$

Das sollte gleich sein wie\\
$Y_{1,1,1}(2. Operator) = -144\,{\frac {1}{ \left( 4\,z_{{1}}+1 \right)  
	\left(192\,z_{{1}}-1 \right) }}$
\\Daraus folgt, dass\\
$\frac{l_x}{l_z} = \left( 4\,x_{{1}}+1 \right)^2$
\\gelten muss, was ich auch schreiben kann als\\
$\frac{l_x}{l_z} = \frac{\left( 4\,x_{{1}}+1 \right)}{\left( - 
	4\,z_{{1}}+1 \right)}$

d.h. als Schnitt von ${\cal L}^2$ transformiert $Y_{1,1,1}(x_{{1}}) mit l_x 
= 1 + 4\,x_{{1}}$ oder $Y_{1,1,1}(z_{{1}})$ mit $l_z = 1 - 4\,z_{{1}}$

Nun kann ich mir auch anschauen, wie sich die holomorphe 3-Form $\Omega$ 
bzw. die Fundamentalperiode $\varpi_0$ des 1. Operators unter dieser 
Variablentransformation verhaelt. Man findet, dass

$\varpi_0(2. Operator) = \varpi_0( z_{{1}} )^2 \left( - 4\,z_{{1}}+1 
	\right) |_{z_{{1}}=\frac {x_{{1}}}{1+4\,x_{{1}}}}$

Das wuerde heissen, dass $\Omega mit \left( - 4\,z_{{1}}+1 \right)^{1/2} $
und ein Schnitt in ${\cal L}^2$ mit $\left( - 4\,z_{{1}}+1 \right) $
transformiert. Das ist aber nicht dasselbe Transformationsverhalten wie 
das von $Y_{1,1,1}$, welches auch ein Schnitt in ${\cal L}^2$ ist. Das 
verstehe ich noch nicht ganz. Vielleicht muss ich da noch ein paar 
weitere Beispiele anschauen.


\bigskip\hrule\bigskip	\let\6=\partial
\[	F(z)=\frac12(\varpi^{(0)}\varpi^{(3)}
		-\sum_1^{h_{21}}\varpi^{(1)}_i\varpi^{(2)}_i)		
\]
\[	Y_{ijk}=\6_{z_i}\6_{z_j}\6_{z_k}F(z)	
\]
Mirror map: $t_i=\varpi^{(1)}_i/\varpi^{(0)}$
\[
	C_{abc}(t)=\frac1{(\varpi^{(0)})^2}\frac{\6z_i}{\6t_a}\frac{\6z_j}
	{\6t_b}\frac{\6z_k}{\6t_c}Y_{ijk}
	=\sum_{d>0} d_ad_bd_cn_d\hbox{Li}_3 q^d
\]
$GW_\beta =\int_{\bar {\cal M}_g(X,\beta)^{vir}}\in\mathbb Q$, 
Instanton numbers = Gopakumar-Vafa = $n_d=\sum_\beta C_g(d,\beta)GW_\beta$
\EndRTransNotes


\section{Outlook}

We have constructed a surprisingly rich set of new Calabi--Yau manifolds 
using conifold transitions from toric Calabi--Yau hypersurfaces.
Small resolutions dual to the flat deformations of the conifold
singularities have been used to construct the mirror families and to
compute quantum cohomologies via mirror symmetry.

The Picard--Fuchs operators have been determined for 28 of the 68 different
diffeomorphism types of one-parameter families. In addition to the 
(computationally expensive)
completion of this calculation it will be interesting to also enumerate
the diffeomorphism types for the large number of cases with 
$h_{11}>1$ and to work out respective Picard--Fuchs operators.
For this, generalizations of the combinatorial formulas for the intersection 
rings and for integral cohomologies should to be derived.
It would also be interesting to extend the calculations of higher-genus 
topological string of \cite{K} to our new one-parameter families and to
check integrality of BPS states as a test for our proposed
mirror construction.

Since already the case of toric hypersurfaces turned out to be a rich
source of new Calabi--Yau 3-folds it would also be interesting to 
generalize our construction to complete intersections. Such transitions,
however of a different type, have already be studied in \cite{Rolling},
where a specialization of the quintic equation and a blow-up of the 
resulting conifold singularity was used to arrive at the
two-parameter bi-degree $(4,1)(1,1)$ complete intersection in 
$\IP^4\ex\IP^1$. In this construction the conifold transition relates
complete intersections of different codimension, but one stays in the
realm of toric ambient spaces.%
\footnote{~
	Note that the topologies of the Calabi--Yau manifolds with 
	$h_{12}=101$ and $h_{12}=103$ in {\it table~\ref{TabPic1}} differ 
	from the toric hypersurfaces with the same Hodge numbers. They 
	originate from polytopes with 8 and 6 vertices and with Picard 
	numbers 4 and 3, respectively.
}

	Another interesting example of codimension two has been discussed in 
appendix E.2 of~\cite{kkrs}, where a toric realization of the hypergeometric 
function 
related to degrees (2,12) with weights (1,1,1,1,4,6), as derived in
\cite{Doran}, is found along a
singular one-parameter subspace of the complex structure moduli
space of a toric complete intersection. 
It is tempting to speculate that a flat deformation of that 
singularity may exist, which could define a smooth Calabi-Yau family with 
the desired monodromy.

{\bf Acknowledgements.} This work is supported in part by the Austrian
Research Funds FWF under grant number P18679-N16. We thank Duco van Straten
and Gert Almkvist for discussions and Emanuel Scheidegger for help in
the computation of instanton numbers.


\end{document}